\documentclass[11pt, letterpaper, oneside]{amsart}

\usepackage{latexsym, amsmath, amssymb, amsfonts, amscd}
\usepackage{amsthm}
\usepackage[mathscr]{eucal}
\usepackage{indentfirst}
\usepackage{enumerate,psfrag}
\usepackage[all,poly,web,knot]{xy}
\usepackage{epsfig, diagrams}

\theoremstyle{plain}
\newtheorem{thm}{Theorem}[section]

\newtheorem{prop}[thm]{Proposition}
\newtheorem{lemma}[thm]{Lemma}
\newtheorem{cor}[thm]{Corollary}

\renewcommand{\latticebody}{\drop@{ }}

\theoremstyle{definition}
\newtheorem{defi}[thm]{Definition}
\newtheorem{pdef}[thm]{Proposition-Definition}

\theoremstyle{remark}
\newtheorem{remark}[thm]{Remark}

\newcommand{\rra}{\Rightarrow}%{\rightrightarrows}

\newcommand{\thra}{\twoheadrightarrow}
\newcommand{\me}{\stackrel{\sim}{\leftrightarrow}} %morita equivalence
\newcommand{\lhc}{\stackrel{\sim}{\to}} %left hypercover
\newcommand{\N}{\ensuremath{\mathbb N}}

\newcommand{\R}{\ensuremath{\mathbb R}}

\newcommand{\D}{\mathscr D}%means differential

\newcommand{\cG}{\mathcal{G}}

\DeclareMathOperator{\colim}{colim}

\newcommand{\bareta}{\bar{\eta}}

\def\U{{\mathcal U}}
\def\L{\Lambda}
\def\D{\Delta}
\newarrow{into} C--->
\newarrow{Eq} =====
\newarrow{dashto} ....>

\def\pD{\partial\D}

\def\PB(#1,#2,#3,#4){
\left\{\begin{matrix}#1&\!\!\!\stackrel{?}{\longrightarrow}&\!\!\!#2\\
\downarrow&&\!\!\!\downarrow\\
#3&\!\!\!\stackrel{?}{\longrightarrow}&\!\!\!#4\end{matrix}\right\}}

\def\pb(#1,#2,#3,#4){ \hom(#1 \to #3, #2 \to #4)}

\def\codi(#1,#2,#3,#4){
\left\{\begin{matrix}#1&\!\!\!\longrightarrow&\!\!\!#2\\
\downarrow&&\!\!\!\downarrow\\
#3&\!\!\!\longrightarrow&\!\!\!#4\end{matrix}\right\}}

\begin{document}

\title{Kan replacement of simplicial manifolds}
\author{Chenchang Zhu \\ Courant Research Centre ``Higher Order Structures'',
University of G\"ottingen}\thanks{Supported by the German Research Foundation 
(Deutsche Forschungsgemeinschaft (DFG)) through 
the Institutional Strategy of the University of G\"ottingen} 
\date{\today}
\maketitle
\begin{abstract}
We establish a functor $Kan$ from local Kan simplicial manifolds to
weak Kan  simplicial manifolds. It gives a solution to the problem of  extending
 local Lie groupoids to Lie $2$-groupoids. 
\end{abstract}

\section{Introduction}

It is a classical topic to study the correspondence between global and infinitesimal 
symmetries. For us, the process from global symmetries to
infinitesimal ones is called differentiation, and the inverse process
is called integration. A classical example of such is in
the case of Lie groups and Lie algebras,
\[ \xymatrix{
& \fbox{\parbox{.3\linewidth}{\center{Lie algebras}}}
  \ar[rrrr]^{\text{differentiation}} &  & & &
  \fbox{\parbox{.3\linewidth}{\center{Lie groups}}}
  \ar[llll]^{\text{integration}} }
\]

However, when our symmetries become more complicated, such as
$L_\infty$-algebras, or even $L_\infty$-algebroids,
the integration and differentiation both become harder. The following
problems have been solved for these higher symmetries: integration of
nilpotent $L_\infty$-algebras by Getzler \cite{getzler}, integration of
general $L_\infty$-algebras by Henriques \cite{henriques}, differentiation of
$L_\infty$-groupoids by \v{S}evera \cite{severa:diff}, both directions for Lie
1-algebroids by Cattaneo-Felder \cite{cafe}, Crainic-Fernandes \cite{cf}, and  from a higher viewpoint by Tseng-Zhu \cite{tz}. Here the author
wants to emphasis  a middle step of local symmetries
missing in the above correspondence,
\[ \xymatrix{
& \fbox{\parbox{.17\linewidth}{\center{Lie algebras}}} \ar[rr]^{\text{local
    integration}}  && \fbox{\parbox{.17\linewidth}{\center{local Lie groups}}}
    \ar[rr]^{\text{extension?}}
  \ar[ll]^{\text{differentiation}} && 
  \fbox{\parbox{.17\linewidth}{\center{Lie groups}}}
  \ar[ll]^{\text{restriction}} }.
\]

Indeed, to do differentiation to obtain infinitesimal symmetries, we
only need local symmetries. Conversely, sometimes, it is easier to obtain a
local integration, avoiding some analytic
issues (for example in \cite{getzler} for $L_\infty$-algebras). In this
paper, we make our first attempts towards the extension
problem from local symmetries to global ones: we construct an extension
from  local Kan simplicial manifolds to weak Kan ones. The classical extension
of local Lie group
 to a topological group discussed by van Est in \cite{van-est:local} can be viewed as 1-truncation
 of our result.
 Its 2-truncation applied to local
Lie groupoids provide a solution to the integration problem of Lie
algebroids
 to Lie 2-groupoids \cite{z:lie2}. Notice that unlike Lie algebras which one-to-one correspond to
simply connected Lie groups, Lie algebroids (integrable or not)
one-to-one correspond to a sort of Lie 2-groupoids with some \'etale property. 

We use the viewpoint of Kan simplicial manifolds to describe arbitrary Lie
$n$-groupoids.

Recall that a simplicial  manifold $X$ consists of  manifolds $X_n$ and structure maps
\begin{equation}\label{eq:fd} d^n_i: X_n \to X_{n-1} \;\text{(face
    maps)}\quad s^n_i: X_n \to X_{n+1} \; \text{(degeneracy maps)},
\;\; \text{for}\;i\in \{0, 1, 2,\dots, n\} 
\end{equation}
that satisfy suitable  coherence conditions (see for example
\cite{friedlander}). The first two examples of simplicial manifolds
(actually, they are simplicial sets with discrete topology) are the simplicial $m$-simplex $\D[m]$ and
the horn $\L[m,j]$ with
\begin{equation}\label{eq:simplex-horn}
\begin{split}
(\D[m])_n & = \{ f: (0,1,\dots,n) \to (0,1,\dots, m)| f(i)\leq
f(j),
\forall i \leq j\}, \\
(\L[m,j])_n & = \{ f\in (\D[m])_n| \{0,\dots,j-1,j+1,\dots,m\}
\nsubseteq \{ f(0),\dots, f(n)\} \}.
\end{split}
\end{equation}
The horn $\L[m,j]$ should be thought as a simplicial set obtained from
$\D[m]$  by taking away its unique
non-degenerate $m$-simplex as well as the $j$-th of its $m+1$
non-degenerate $(m-1)$-simplices. 

\vspace{.6cm}
\centerline{\epsfig{file=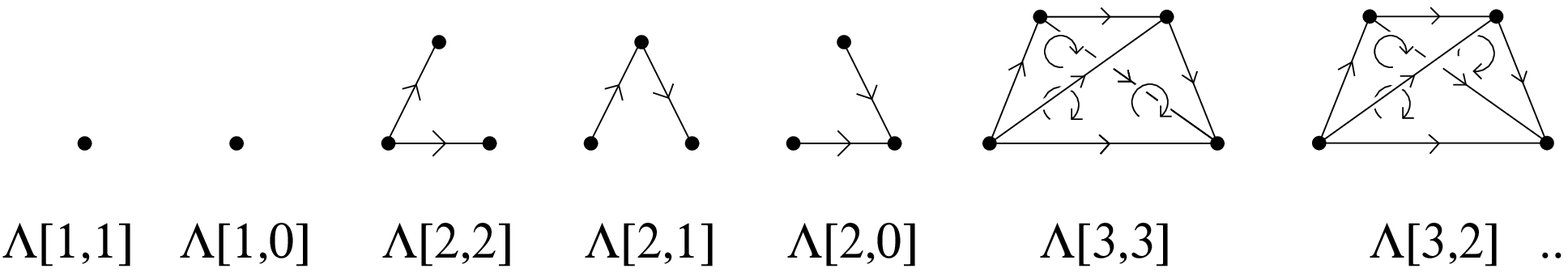,height=1.9cm}}
%\vspace{.3cm}
Our convention for arrows is that they are oriented from bigger numbers to
smaller numbers.

Let us recall that in homotopy theory, Kan conditions say that the
natural restriction map 
\begin{equation}\label{eq:horn-proj}
X_m =\hom(\D[m], X) \to \hom(\L[m,j], X).
\end{equation}
is surjective, i.e. any horn can be filled up by a simplex. They correspond to the
possibility of composing and inverting various morphisms, in the
language of groupoids. 

With enrichment in differential geometry,  {\bf Kan conditions} are 

\vspace{.2cm}
\noindent
\begin{tabular}{p{1.8cm}p{5.8cm}p{1.8cm}p{5.8cm}}
$Kan(m,j)$:& \eqref{eq:horn-proj}  is a
surjective submersion, &
$Kan!(m,j)$:& \eqref{eq:horn-proj} is a
diffeomorphism.\\
\end{tabular} 
\vspace{.2cm}

But since $\hom(\L[m,j], X)$ is formed by taking a numerous fibre products
of the $X_i$'s, it may not be a manifold. However 
if \eqref{eq:horn-proj}
is a submersion for all $0\le j \le m \le m_0$, as shown in
\cite[Lemma 2.4]{henriques}, $\hom(\L[m_0,j], X)$
is a manifold for all $0\le j \le m_0$. Hence we are allowed to
define,

\begin{defi}\label{def:defngroupoid} A {\em Lie} $n${\em-groupoid} $X$
($n\in\N \cup \infty$) is a simplicial manifold that satisfies
$Kan(m,j)$, $\forall m\ge 1$, $0\le j\le m$, and $Kan!(m,j)$
,$\forall m> n$, $0\le j\le m$.
When $n=\infty$, a Lie $\infty$-groupoid is also called a {\em Kan
  simplicial manifold}. 
\end{defi}
Then Lie $1$-groupoid is simply the nerve of a Lie groupoid.

To describe local
Lie groupoids, we need  {\bf local Kan conditions}:

\begin{tabular}{p{1.8cm}p{3.8cm}p{1.8cm}p{3.8cm}}
$Kan^l (m,j)$:& \eqref{eq:horn-proj}  is a submersion, &
$Kan^l !(m,j)$:& \eqref{eq:horn-proj} is injective \'etale.\\
\end{tabular} 
%\vspace{.2cm}

\begin{defi}\label{def:local-ngroupoid} A local {\em Lie} $n${\em-groupoid} $X$
($n\in\N \cup \infty$) is a simplicial manifold that satisfies
$Kan^l(m,j)$, $\forall m\ge 1$, $0\le j\le m$, and $Kan^l!(m,j)$
,$\forall m> n$, $0\le j\le m$. When $n=\infty$, a local Lie
 $\infty$-groupoid is also called a  {\em local Kan
  simplicial manifold}. 
\end{defi}

Then a local Lie $1$-groupoid $X$ is the nerve of a local Lie
groupoid. 

As soon as we have done this, it becomes clear that to associate a
Kan object $Kan(X)$ to a local Kan simplicial manifold $X$, we need to do some sort of
fibrant replacement in the category of simplicial manifolds. However,
simplicial manifolds do not form a model category and we need to
do it by hand. In fact, the differential category is rather special, even the
construction for simplicial presheaves can not be used directly
here. It turns out that the object $Kan(X)$ constructed directly by Quillen's small
object argument is not a Kan simplicial manifold, however it is a
simplicial manifold and is Kan as a simplicial set. We also prove certain
representibility conditions for $Kan(X)$ and make it into a weak Kan
simplicial manifold (see Section \ref{sec:defi}), which is slightly
weaker than a Kan one. On the other hand, the defects of $Kan(X)$ lie only on high levels,
that is, if we perform a 2-truncation $\tau_2(Kan(X))$, and the
2-truncation  is still representable, then   $\tau_2(Kan(X))$ is
indeed a Lie 2-groupoid.

\section{Definition} \label{sec:defi}

Now we try to define a functor $Kan$ sending invertible local Kan
manifold to Kan simplicial manifolds by modifying directly Quillen's small object
argument. We will see that it is  not successful, however we
arrive at a simplicial manifold satisfying conditions slightly weaker
than Kan.   Let
 \begin{equation}\label{eq:j} J:= \{ \Lambda[k, j] \to \Delta[k]: 0\le j \le k\ge 3,  \} \cup \{ \Lambda[2, 1]
\to \Delta[2] \}, \end{equation}  
be a subset
of inclusions with respect to which Kan condition have the right lifting property.
Given a local Kan manifold $ X$, we then construct a series of simplicial manifolds
\begin{equation}\label{eq:series}
X=X^0 \to X^1 \to X^2 \to \dots \to X^\beta \to \dots  
\end{equation}
by an inductive push-out:
\begin{equation}\label{diag:x}
\begin{diagram}
\coprod_{(\Lambda[k, j] \to \Delta[k]) \in J} \Lambda[k, j] \times
\hom(\Lambda[k, j], X^\beta)   &    \rTo    &  X^\beta \\
 \dTo       &   & \dTo \\
\coprod_{(\Lambda[k, j] \to \Delta[k]) \in J} \Delta[k]\times \hom(\Lambda[k, j], X^\beta)
  &    \rTo    &  {\NWpbk} X^{\beta+1}&.
\end{diagram}
\end{equation}
Then we let $ Kan(X)= \colim_{\beta \in \N} X^\beta$.

Now we make some calculation for first several steps of Kan replacement:
First of all $X_0 = X^1_0 = X^2_0 = \dots
= Kan(X)_0$, and 
\begin{equation}\label{eq:x1}
\begin{split}
X^1_1 &=  X_1 \sqcup (X_1 \times_{X_0} X_1) \\
X^2_1 &= X^1_1 \sqcup X^1_1 \times_{X_0} X^1_1 \\
      &= X_1^1 \sqcup \big( X_1 \times_{X_0} X_1
\sqcup X_1 \times_{X_0} (X_1 \times_{X_0} X_1) \\
&\sqcup (X_1 \times_{X_0} X_1)  \times_{X_0} X_1 \sqcup  (X_1
\times_{X_0} X_1)\times_{X_0}(X_1 \times_{X_0} X_1) \big) \\
& \vdots \\
 Kan(X)_1 &= X_1 \sqcup ( X_1 \times_{X_0} X_1) \sqcup (X^1_1
 \times_{X_0} X^1_1 ) \sqcup  (X^2_1 \times_{X_0} X^2_1)  \dots,  
\end{split}
\end{equation}
which we can represent them by the following picture: \\
\psfrag{...}{$\dots$} \psfrag{X1}{$Kan(X)_1:$}
\centerline{\epsfig{file=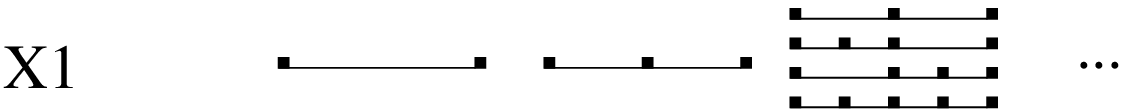, height=1cm, width=8cm}}

A calculation shows that
\[
\begin{split}
X^1_2 &= X_2 \sqcup X_1 \times_{X_0} X_1 \sqcup X_1 \times_{X_0} X_1
\sqcup X_1 \times_{X_0} X_1 \\
& \sqcup (\sqcup_{j=0}^3\hom(\L[3, j], X) \\
X^2_2 &= X^1_2 \sqcup X^1_1 \times_{X_0} X^1_1 \sqcup X^1_1 \times_{X_0} X^1_1
\sqcup X^1_1 \times_{X_0} X^1_1 \\
& \sqcup (\sqcup_{j=0}^3\hom(\L[3, j], X^1) \\ 
&\vdots
\end{split}
\]
Inside $X^1_2$, there are three copies of $X_1\times_{X_0} X_1$. The
first is an artificial filling of the horn $X_1\times_{X_0} X_1$, and
the second two are images of degeneracies of $X_1\times_{X_0} X_1$ in
$X^1_1$. The same for $X^2_2$, etc.  We represent an element in $X^1_2$
as\\
\psfrag{X2}{$X_2:$}\psfrag{X11}{$X_1
  \times_{X_0}X_1:$}\psfrag{hom}{$\hom(\L(3, j), X):$}
\psfrag{...4 of them}{$\dots$ 4 such} \psfrag{,}{,}
\begin{equation}\label{pic:x12}
\epsfig{file=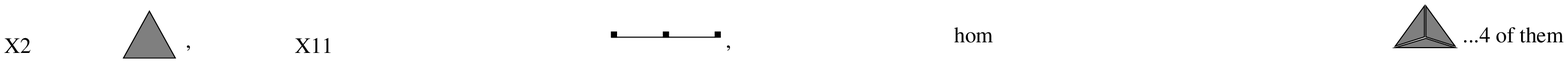, height=.5cm, width=12cm }
\end{equation}
plus those degenerate ones in the other two copies of  $X_1
\times_{X_0} X_1$. Furthermore we represent an element
in $X^2_2$ as 

\vspace{.5cm}
\psfrag{X21}{$X_2^1$: described as above}
\psfrag{X11}{$X_1^1\times_{X_0}X_1^1$:} \psfrag{X3j}{$\hom(\L[3, j],
  X^1)$:} \psfrag{... 7 of them}{$\dots$ }
\begin{equation}\label{pic:x22}
\epsfig{file=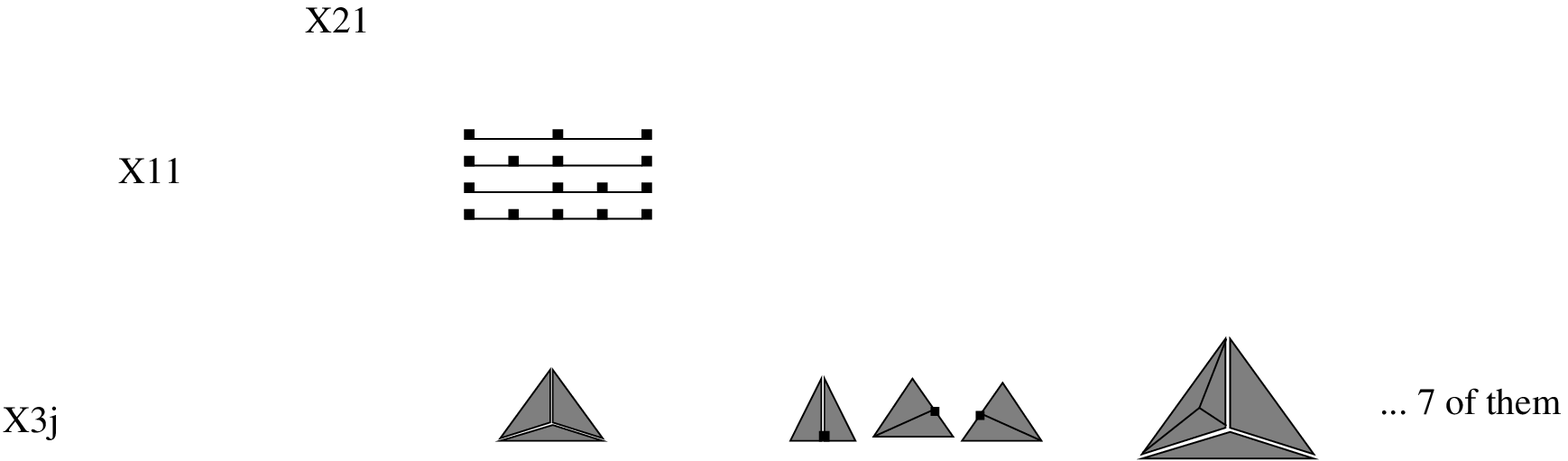, height=3cm, width=12cm}
\end{equation}

We now show that
 \begin{equation}\label{eq:x2-to-horn}
X_2^1 \to \hom(\Lambda[2, 0], X^1)=X^1_1 \times_{d_2,
  X_0, d_1}
X^1_1\end{equation} 
is actually not a submersion. 

We first need some technical preparation. A simplicial set $S$ is {\em collapsible} if it admits
a filtration 
\begin{equation}\label{eq:collapsible}
 pt=S_0 \subset S_1 \subset \dots \subset S_k= S \end{equation} 
such that each $S_i$ is obtained from the previous one by filling a
horn, namely such that $S_i$ can be written as $S_i = S_{i-1}
\sqcup_{\L[n_i, l_i]} \D[n_i]$ for some injective map $\L[n_i, l_i] \hookrightarrow
S_{i-1}$.  Thus we have an order for collapsible simplicial
sets: we say $S$ is not bigger than $T$, denoted as $S\prec T$,  if
$T=S_t $ and $S=S_s$  with $s\leq t$ in
\eqref{eq:collapsible}. For us the notation $S\prec T$ also indicates the
inclusion map $S\to T$. We also define the dimension of a collapsible
simplicial set $S$ as 
\[ \dim S = \max_k\{ \Delta[k] \prec S \}.  \]

\begin{lemma}\label{lemma:st}
Given $S \prec T$ and a local Kan
simplicial manifold $X$, $\hom(T,X)$ and $\hom(S, X)$ are both
manifolds, and the natural map
\[ \hom(T, X) \to \hom(S, X), \]
is always a submersion.
\end{lemma}
This is proven in \cite[Lemma 2.4]{henriques}
for Kan simplicial manifolds, but it is easy to see that it works also
for local Kan ones since only the submersion condition is used.

Back to the map \eqref{eq:x2-to-horn}, $X_2^1$ has several components, and the horn
projection map \eqref{eq:x2-to-horn} induces on each component the
following maps:
\begin{eqnarray} 
& X_2 \to \hom(\Lambda[2,0], X), \qquad  & \hom(\Lambda[3,j], X) \to
\hom(\Lambda[2,0], X)  \label{eq:non-deg}\\
& X_1\times_{X_0} X_1 \to X_1 \times_{X_0} (X_1 \times_{X_0} X_1), 
& X_1\times_{X_0} X_1 \to (X_1 \times_{X_0} X_1) \times_{X_0} (X_1
\times_{X_0} X_1), \label{eq:deg1}\\
& X_1\times_{X_0} X_1 \to X_1 \times_{X_0} X_1 \label{eq:deg2}.    
\end{eqnarray}
The morphisms in \eqref{eq:non-deg} are always submersions by Lemma
\ref{lemma:st}. However the morphisms in \eqref{eq:deg1} are not
submersions. This implies that $\hom(\D[2],
Kan(X)) \to \hom(\Lambda[2,0], Kan(X))$ will not be a
submersion. Hence $Kan(X)$ will not be a Kan simplicial manifold, 
but we will prove that it satisfies
\begin{enumerate}
\item [(A)] $Kan(X)$ is a simplicial manifold;
\item [(B)] moreover, $\hom(S, Kan(X))$ is a manifold for any collapsible $S$; 
\item [(C)] the natural map $Kan(X)_n \to \hom(\Lambda[n, l], Kan(X))$ is
  surjective for all $n$ and $l$ with $0\leq l \leq n$. 
\end{enumerate}
Simplicial manifolds satisfying such conditions are called {\bf weak
  Kan simplicial manifolds}. These weak Kan simplicial
manifolds as simplicial sets are indeed Kan. The submersion condition
in Kan condition
is replaced by condition (B). By Lemma \ref{lemma:st}, we can see that
the submersion condition implies (B), hence weak Kan is indeed weaker
than Kan. However, the usual place to use submersion condition is to
guarantee  some representibility for example the one in condition
(B). Hence we see that in many cases we can bypass the difficulty thanks
to this condition. For example, we can still talk about hypercovers of
these weak Kan simplicial manifolds (even though it is not used in
this paper).

Before attacking the problem, we first prepare a technical lemma:

\begin{lemma}\label{lemma:1-0}
Suppose that $X^\beta$'s are a sequence of simplicial sets constructed by 
\eqref{diag:x},
\begin{enumerate}
\item \label{itm:s} if $S$ is a collapsible simplicial set, then we
  can decompose 
\[\hom(S, X^{\beta+1}) = \sqcup_{a\in
    A} \hom(S_a , X^\beta ),\] with a finite set of  collapsible
  simplicial sets $\{S_a: a \in A\}\ni S$ satisfying $ \dim S_a \leq \dim S $. 
\item \label{itm:st} this decomposition respects morphisms, that is if
  $S\prec T$, and both of them have a
  decomposition, 
\[\hom(S, X^{\beta+1}) = \sqcup_{a\in
    A} \hom(S_a , X^\beta ), \quad \hom(T, X^{\beta+1}) = \sqcup_{a'\in
    A'} \hom(T_{a'} , X^\beta ),
\] 
 then there is a map $a: A' \to A$, and  morphisms of simplicial sets
 $S_{a(a')} \to T_{a'} $, such that the natural morphism $\hom(T,
 X^{\beta+1}) \to \hom(S,
 X^{\beta+1})$ is induced from $\hom(T_{a'} , X^\beta) \to \hom
 (S_{a(a')}, X^\beta)$ on the level of their decompositions.
\end{enumerate}
\end{lemma}
\begin{proof}
Since the procedure to form $X^{\beta+1}$ by $X^\beta$ is the same
as the one to form $X^1$ by $X$,  we only have to prove the two statements for $\beta=0$. Since $X$ is arbitrary, the decomposition in \eqref{itm:s}
is clearly unique. We use an
induction on the size of $S$ and
$T$. The initial assumption is  verified in the calculation we did
earlier in this section. It is
clear that \eqref{itm:s} holds for $\D[m]$ for $m\in \N$.  Now with a fixed
$n$, we consider a horn filling diagram as we mentioned in the process
of \eqref{eq:collapsible} 
\begin{equation}\label{eq:st} 
\begin{diagram}
S   &    \rTo    &  {\SWpbk} T \\
 \uTo       &   & \uTo \\
\Lambda[k,j]  &    \rTo    &   \Delta[k]&,
\end{diagram}
\end{equation}
with $\dim S \leq \dim T \leq n-1$, and $k\leq n-1$. 

We suppose that
\begin{enumerate}
\item[(i)] \label{itm:ss'} statement \eqref{itm:s} is true for all
  $S'$ with $S' \prec S$;
\item[(ii)] \label{itm:s't'} statement \eqref{itm:st} is true for $S' \prec T'$ and $S' \prec
  \D[n]$ when $\dim S' \leq \dim T' \leq n-1$ and when \eqref{itm:s}
  verifies for $S'$ and $T'$.
\end{enumerate}
To finish the induction, we will prove that
\begin{itemize}
\item statement \eqref{itm:s} holds for $T$;
\item statement \eqref{itm:st} holds for $S\prec T$, and $T\prec
  \Delta[n]$ if such a map $T\to \D[n]$ exists. 
\end{itemize}

First of all, we apply $\hom(-, X^1)$ to \eqref{eq:st} and apply the
induction hypothesis to $\hom(S, X^1)$, $\hom(\Lambda[k,j], X^1)$ and
$\hom(\D[k], X^1)$, then we have
\[\begin{split} 
\hom(T, X^1) &= \hom(S, X^1) \times_{\hom(\Lambda[k,j], X^1)}
\hom(\D[k], X^1) \\
&= \sqcup_{a'\in A'} \hom(S'_{a'} , X) \times_{\sqcup_{a\in A}
  \hom(S_a, X)} \sqcup_{a''\in A''} \hom(S''_{a''}, X) \\
&=\sqcup_{b\in B} \hom(T_b, X).
\end{split}
\]
Here $T_b$ is formed when $a(a')=a(a'')$ by
\[
\begin{diagram}
S'_{a'}   &    \rTo    &  {\SWpbk} T_b \\
 \uTo       &   & \uTo \\
S_{a(a')}  &    \rTo    &   S''_{a''}&.
\end{diagram}
\]
We obtain a map $B\to A'$ defined by $b\mapsto a'$ and morphisms $S'_{a'}
 \to T_b$. They induce the morphisms $\hom(T_b, X) \to
 \hom(S'_{a'}, X)$, hence the morphism
 $\hom(T, X^1) \to \hom(S, X^1)$.  It's not hard to see that $T\in \{T_b\}$
 by induction hypothesis and \eqref{eq:st}. 

Suppose
 $\hom(\D[n], X^1) =\sqcup_{c\in C} (D_c, X) $. If there is a map
 $T\prec \D[n]$, by restriction, we obtain maps $S \prec \D[n]$,
 $\Lambda[k, j]\prec \D[n]$, and $\D[k] \prec \D[n]$ which fit in the
 following commutative diagram:
\[
\begin{diagram}
S   &    \rTo    &   \D[n] \\
 \uTo       & \ruTo   & \uTo \\
\Lambda[k,j]  &    \rTo    &   \D[k]&,
\end{diagram}
\]
By induction hypothesis, we have
\begin{itemize}
\item the morphism $\hom(\D[n], X^1) \to \hom(S, X^1)$ is 
 induced by a map $a': C\to A'$ and morphisms $S'_{a'(c)} \to D_{c}$;
\item   the morphism $\hom(\D[n], X^1) \to \hom(\Lambda[k,j], X^1)$ is
 induced by a map $a: C \to A$ and morphisms $S_{a(c)} \to D_{c}$;
 \item  the morphism $\hom(\D[n], X^1) \to \hom(\Delta[k], X^1)$ is
 induced by a map $a'': C \to A''$ and  morphisms $S''_{a''(c)} \to D_{c}$.
\end{itemize} 
We see that $ \hom(D_c, X) \to \hom(S_{a(c)}, X)$  induces $\hom(\D[n], X^1) \to \hom(\Lambda[k, j], X^1)$, 
and the composed morphism $\hom(D_c, X) \to \hom(S'_{a'(c)}, X) \to
\hom(S_{a(a'(c))}, X)$  induces $ \hom(\D[n], X^1) \to \hom(S, X^1) \to
\hom(\Lambda[k, j], X^1)$, 
which is the same morphism as $\hom(\D[n], X^1) \to \hom(\Lambda[k,
j], X^1)$. Hence by uniqueness of the decomposition, we have $a(a'(c))=a(c)$ and similarly $a(c)=a(a''(c))$,
and a commutative diagram
\[
\begin{diagram}
S'_{a'(c)}   &    \rTo    &   D_c \\
 \uTo       & \ruTo  
 & \uTo \\
S_{a(a')}  &    \rTo    &   S''_{a''(c)}&,
\end{diagram}
\]

Then $T_{b(c)}$ defined by the
pushout diagram 
\[
\begin{diagram}
S''_{a(a''(c))}   &    \rTo    &  {\SWpbk} T_{b(c)} \\
 \uTo       &   & \uTo \\
S_{a(c)}  &    \rTo    &   S'_{a(a'(c))}&.
\end{diagram}
\] 
has a canonical map $T_{b(c)} \to D_c$. By the property of $\hom(T,
X^1)$ being the fibre
product,  these canonical maps induce the map
$\hom(\D[n], X^1) \to \hom(T, X^1) $ via the maps $\hom(D_c, X) \to
\hom(T_{b(c)}, X) $. 
\end{proof}

\begin{pdef} \label{pdef:kani} The operation $Kan$ constructed in
 \eqref{diag:x} is a
functor  from the category of local Kan manifolds $
  X$  to the one of weak Kan simplicial manifolds.
\end{pdef}
\begin{proof}
The construction
of $Kan$ makes it clear that it  is
functorial. Since $Kan(X) =
\colim_{\beta} X^\beta$, given any finite simplicial set $A$ (a collapsible
simplicial set $S$ is such), the natural map of sets is an
isomorphism,
\begin{equation}\label{eq:iso-colim}\colim_{\beta} \hom(A, X^\beta)
  \xrightarrow{\simeq} \hom(A, Kan(X)). 
\end{equation}
 Moreover by Lemma \ref{lemma:1-0}, 
\begin{equation}\label{eq:xn}
\hom(S, X^{\beta+1})=  \hom(S, X^\beta)  \bigsqcup (\sqcup_a \hom(S_a , X^\beta)),
\end{equation}
We then use Lemma \ref{lemma:1-0} recursively, and obtain that for any
collapsible simplicial set $T$,  
\[ \hom(T, X^{\beta})= \sqcup \hom(T_p, X), \] for a finite set of
collapsible simplicial sets $T_p$. Hence $\hom(S,
X^{\beta})$ and $\hom(S_a , X^\beta)$ are manifolds because $X$ is
local Kan. By \eqref{eq:xn} and \eqref{eq:iso-colim},  $\hom(S, Kan(X))$ is a disjoint union of manifolds.

So it remains to show that $Kan(X)$ is Kan as a simplicial set.  We take an element $A\to B$ of $J$ and a solid arrow
diagram,
\begin{equation}
\label{eq:surj}
\xymatrix{
A \ar[r] \ar[d] & Kan(X)  \ar[d] \\
B \ar[r] \ar@{.>}[ru]  & pt
}
\end{equation}
then we must show that the dotted arrow exists. By the isomorphism \eqref{eq:iso-colim}, the map $A\to Kan(X)$ factors through $X^\beta \to
Kan(X)$ for some $\beta $ and we have the solid arrow diagram
\[
\xymatrix{ A \ar[r] \ar[d] & X^\beta \ar[r] \ar[d] & X^{\beta+1}
  \ar[r] \ar[ld] & Kan(X) \ar[lld] \\
B \ar@{.>}[rru] \ar[r] & pt
} 
\]
Since $X^{\beta+1}$ is constructed as the push-out in \eqref{diag:x},
the dotted arrow naturally exists, and this dotted arrow defines the
one in \eqref{eq:surj}.

Now we only have to verify that the dotted arrow in \eqref{eq:surj} exists for
$\L[1,j]\to \D[1]$ for $j=0,1$ and $\L[2, j] \to \D[2]$ for $j=0,2$.
We have $X^\beta_0 = X_0$, and
\[\hom(\L[1, j], X^{\beta+1})= X^{\beta+1}_0 = X_0, \quad \hom (\D[1],
X^{\beta+1}) = X^\beta_1 \times_{d_0, X_0, d_1} X^\beta_1, \]
thus the map $\hom (\D[1],
X^{\beta+1})\to \hom(\L[1, j], X^{\beta+1}) $ being the pull-back of
$d_1$ or $d_0$, has to be a surjective submersion. 
Now we prove that if $X^\beta$ is invertible, then $X^{\beta+1}$ is
also invertible. 
\[
\begin{split}
 &\hom (\Lambda[2,2] , X^{\beta +1})\\ = &X^{\beta+1}_1 \times_{d_1, X_0,
  d_1} X^{\beta+1}_1  \\  =&( X^{\beta}_1 \sqcup  X^{\beta}_1\times_{d_0,
  X_0, d_1}  X^{\beta}_1) \times_{d_1, X_0, d_1} ( X^{\beta}_1 \sqcup  X^{\beta}_1\times_{d_0,
  X_0, d_1}  X^{\beta}_1) \\
=& X^{\beta}_1\times_{d_1,
  X_0, d_1}  X^{\beta}_1 \sqcup   X^{\beta}_1\times_{d_1,
  X_0, d_1}  ( X^{\beta}_1\times_{d_0,
  X_0, d_1}  X^{\beta}_1) \sqcup (X^{\beta}_1\times_{d_0,
  X_0, d_1}  X^{\beta}_1)\times_{d_1,
  X_0, d_1} X^{\beta}_1 \sqcup \dots
\end{split} 
\]
Since $X^\beta$ is invertible, $ X^{\beta}_1\times_{d_1,
  X_0, d_1}  X^{\beta}_1\cong X^{\beta}_1\times_{d_0,
  X_0, d_1}  X^{\beta}_1 $. Hence
\[ 
\begin{split}
  X^{\beta}_1\times_{d_1,X_0, d_1}  ( X^{\beta}_1\times_{d_0,X_0, d_1}
  X^{\beta}_1) & \cong X^\beta_1 \times_{d_0, X_0, d_1} X^{\beta}_1\times_{d_0,X_0, d_1}
  X^{\beta}_1, \\
(X^{\beta}_1\times_{d_0,
  X_0, d_1}  X^{\beta}_1)\times_{d_1,
  X_0, d_1} X^{\beta}_1  &\cong X^{\beta}_1\times_{d_1,
  X_0, d_1}  X^{\beta}_1 \times_{d_0,
  X_0, d_1} X^{\beta}_1 \cong (X^{\beta}_1\times_{d_0,
  X_0, d_1}  X^{\beta}_1) \times_{d_0,
  X_0, d_1} X^{\beta}_1, \\ & \dots
\end{split}
\]
It is easy to continue to verify that $X^{\beta+1}$ is invertible. Then the
final result follows from  \eqref{eq:iso-colim}. 
\end{proof}

Given an invertible local Kan manifold $X$, we call $Kan(X)$ the {\bf Kan
  replacement} of $X$. 

Even through $Kan(X)$ is not Kan, its 2-truncation $\tau_2(Kan(X))$
behaves well. We define $n${\em -truncation}
$\tau_n$ (it is called
$\tau_{\le n} $ in \cite[Section 3]{henriques}), 
of a simplicial manifold $X$ as, 
\[  \tau_n(X)_k = X_k, \forall k\le n-1, \quad \tau_n(X)_k =
X_k/\sim_k, \forall k \ge n, \]
where two elements $x\sim_k y $ in $X_k$ if they are
homotopic\footnote{This means that $d_i x = d_i y$, $0\le i \le k$,
  and there exists  $z\in X_{k+1}$ such that $d_k(z)=x, d_{k+1}(z)=y$,
and $d_i z = s_{k-1} d_i x = s_{k-1} d_i y$, $0\le i < k$. } and have the
same $n$-skeleton. Since in the procedure, taking a quotient is
involved, the result $\tau_n(X)$ might not be a simplicial manifold
anymore. We view it as a simplicial stack. When $X$ is Kan,  $\tau_n(X)$
viewed as a simplicial set is always a discrete $n$-groupoid. It is representable, namely
it is indeed a simplicial manifold, if and only if the quotient
$X_n/\sim_{n}$ is representable because the higher levels are decided
by $X_n/\sim_n$.  Even though $Kan(X)$ is not a Kan manifold, we still have

\begin{prop}\label{prop:tau2}
When $Kan(X)_2/\sim_2$ is representable, $\tau_2(Kan(X))$ is a Lie
2-groupoid. 
\end{prop}
\begin{proof}
As a simplicial set, $Kan(X)$ is Kan. Hence $\hom(\D[n],\tau_2(Kan(X))) \cong \hom(\Lambda[n,j],
\tau_2(Kan(X))) $, for $n\geq 3$. Especially, $\tau_2(Kan(X))_3\cong
\hom(\Lambda[3,0], \tau_2(Kan(X)))$. Since the higher layers are
determined by the first four layers,
\[ \tau_2(Kan(X))=Cosk^3 \circ Sk^3(\tau_2(Kan(X))), \]
by the same argument in \cite[Section 2.3]{z:tgpd}, to
show $\tau_2(Kan(X))$ is a Lie 2-groupoid, we only need to show that
$\hom(\Lambda[3,0], \tau_2(Kan(X)))$ is representable and $Kan(m\le 2,j)$
for $\tau_2(Kan(X))$. In fact the induction argument there already
shows that the representibility of
$\hom(\Lambda[3,0], \tau_2(Kan(X)))$  is implied by $Kan(m\le 2,j)$
for $\tau_2(Kan(X))$   given $\tau_2(Kan(X))_2=Kan(X)_2/\sim_2$ is
representable. Hence we only need to show  $Kan(m\le 2, j)$.

As shown in Def.-Prop. \ref{pdef:kani}, $\hom (\D[1],
X^{\beta+1})\to \hom(\L[1, j], X^{\beta+1}) $ being the pull-back of
$d_1$ or $d_0$, is a surjective submersion, hence $\hom(\D[1], Kan(X))
\to \hom(\L[1,j], Kan(X))$ is a surjective submersion. This is $Kan(1,j)$ for $Kan(X)$, hence for $\tau_2(Kan(X))$. 

The surjective part in $Kan(2, j)$ is automatically satisfied: since
$Kan(X)$ is Kan as a simplicial set,  the composed map
\[Kan(X)_2 \to \tau_2(Kan(X)) \xrightarrow{p} \hom(\Lambda[2, j], Kan(X))=
\hom(\Lambda[2, j], \tau_2(Kan(X))), \]
is surjective, hence the desired map $p$ is also surjective. We only need to show the
submersion part.  Then what happened to the degenerate faces  where the
horn projection map
is not a submersion for $Kan(X)_2$? An element  $\eta \in Kan(X)_2$
can be described as a tree as stated in Lemma \ref{lemma:tree}. If all
the vertices of the tree are triangles in $X_2$, then the horn
projection map is a submersion for $Kan(X)_2$, hence $p$ is a submersion. The problem happens exactly when the tree contains
at least one vertex coming from one of the three copies of  $X_1 \times_{X_0} X_1$.  But  these bad pieces as in
\eqref{eq:deg1} and \eqref{eq:deg2} are all homotopic via elements in
$Kan(X)_3$ to the boundary of good
pieces as in \eqref{eq:non-deg}, where the submersion holds.  Hence the submersion part is also
true for the 2-truncation.

\end{proof}

\section{Universal Properties}

Given a local Lie $1$-groupoid $W$ (or the nerve of a local Lie
groupoid), then it extends to a Lie $2$-groupoid 
$\tau_2(Kan(W))$. In \cite{z:lie2}, we verified that
$\tau_2(Kan(W))$ is always a Lie $2$-groupoid (even though
$\tau_1(Kan(W))$ might not be Lie) with universal property.

For this purpose, we need to show some universal properties of our Kan replacement. It should be stable under Morita equivalence of simplicial
manifolds (whatever that is), and if some simplicial manifold $X$ is
already Kan, $Kan(X)$ should be Morita equivalent to $X$. Hence let's
first begin with an introduction of these concepts such as Morita
equivalence. 

\subsection{Morita equivalence of local Kan manifolds} \label{sec:me}

The reader's first guess is probably that a morphism $f:X\to Y$ of
 simplicial manifolds ought to be a simplicial smooth map i.e. a collection of
smooth maps $f_n:X_n\to Y_n$ that commute with faces and
degeneracies. We
shall call such a morphism a \em strict map \rm from
$X$ to $Y$. Unfortunately, it is known that, already in the case
of usual Lie groupoids, such strict notions are not good enough.
Indeed there are strict maps that are not invertible even though
they ought to be isomorphisms. That's why people introduced the
notion of {\em Hilsum-Skandalis bimodules} \cite{hs}. Here is an
example of such a situation: consider a manifold $M$ with an open cover
$\{\U_\alpha\}$. The simplicial manifold $X$ with
$X_n=\bigsqcup_{\alpha_1,\ldots,\alpha_n}\U_{\alpha_1}\cap\cdots\cap\U_{\alpha_n}$
maps naturally to the constant simplicial manifold $M$. All the
fibers of that map are simplices, in particular they are
contractible simplicial sets. Nevertheless, that map has no
inverse.

The second guess is then to define a special class of strict maps
which we shall call {\em hypercovers}. A map from $X$ to $Y$
would then be a {\em zig-zag}  of strict maps
$X\stackrel{\sim}{\leftarrow}Z\to Y$, where the map $Z\to X$ is
one of these hypercovers.

Another alternative however equivalent way to define a generalized morphism of
simplicial 
manifolds follows from \cite[Section 2.4]{lurie}'s Cartesian
fibrations. In this paper, we use the zig-zag method with the notion of
hypercover.

Our hypercover is very much inspired from the notion of
hypercover of \'etale simplicial objects \cite{SGA4, friedlander}  and of
trivial fibration of Quillen for simplicial sets \cite{quillen:ha}.

Recall \cite[Section I.3]{may}, given a pointed Kan simplicial set $X$, i.e. $X_0 = pt$, its homotopy groups are given by 
\[
 \pi_n(X):=\{x\in X_n | d_i(x)=pt  \;\text{for all} \; i\}/\sim
\]
where $x\sim x'$ if there exists an element $y\in X_{n+1}$ such that
$d_0 (y)= x $, $d_1(y)=x'$, and $d_i(y)=pt$ for all $i>1$. When $X_0$ is
not necessarily a point, $\pi_n$ is a sheaf over $X_0$ in general.

\begin{lemma}\label{eqss}%<eqss>
Given a map  $S\to T$ of pointed Kan simplicial sets,  if for
any $n\ge 0$ and any commutative solid arrow diagram
\begin{equation}\label{eqeq}
\xymatrix{
\partial\D[n] \ar[r] \ar@{^{(}->}[d] &S \ar[d]\\
\D[n]\ar[r] \ar@{.>}[ru] & T}
\end{equation}
there exists a dotted arrow that makes both triangles commute, then this map is a homotopy equivalence,
i.e. $\pi_n(S)=\pi_n(T)$.
Here $\partial \D[n] $ stands for the boundary of the $n$-simplex. 
\end{lemma}
The proof is standard.

Translating the condition of Lemma \ref{eqss} into hom spaces
gives:

\begin{defi}\label{defequivalence}
A strict map $f:Z\to X$ of {\em local Kan simplicial manifolds  is a
 hypercover}
 if the natural
map
\begin{equation}\label{eq:pb} 
Z_m=\hom(\D[m],Z)  \to \hom(\pD[m]\to \D[m], Z\to X) \end{equation} 
is a surjectve submersion for all $0\le m$. 
\end{defi}

Here $ \hom(A\to B, Z\to X)$  denotes the pull-back spaces of the form
$\hom(A,Z)\times_{\hom(A,X)}\hom(B,X)$, where the
maps are induced by some fixed maps $A\to B$ and $Z\to X$.
This notation indicates that the space parameterizes all commuting
diagrams of the form
$$\begin{matrix}A&\!\!\!\longrightarrow&\!\!\!Z\\
\downarrow&&\!\!\!\downarrow\\
B&\!\!\!\longrightarrow&\!\!\!X,\end{matrix}$$ where we allow the
horizontal arrows to vary but we fix the vertical ones.

Similarly, we can define hypercover for Lie $n$-groupoids:
\begin{defi}\label{defequivalence-gpd}
A strict map $f:Z\to X$ of {\em Lie $n$-groupoids is a hypercover} 
if the natural
map \eqref{eq:pb} is a  surjective submersion for all $0 \le m< n$ and is an isomorphism
when $m=n$.
\end{defi}
\begin{remark}
As proved in \cite{z:tgpd-2}, if $f: Z\to X$ is a hypercover of Lie
$n$-groupoids, then   \eqref{eq:pb} is
automatically an  isomorphism for all $m>n$. 
\end{remark}

As in the case of Definition \ref{def:defngroupoid}, we need to
justify that the pull-back $ \hom(\pD[m]\to \D[m], Z\to X) $ is a
manifold. 
This is rather surprising since the spaces
$\hom(\pD[m],Z)$ need not be manifolds (for example take $m=2$ and
$Z$ the cross product groupoid associated to the action of $S^1$
on $\R^2$ by rotation around the origin). We
justified this in \cite{z:tgpd-2} for Kan simplicial manifolds, but it is
clear that only the submersion property is needed, hence the same proof
works for local Kan manifolds.

\begin{defi}\label{defi:m-equi-2gpd} Two local Kan simplicial manifolds $X$ and
  $Y$ are
{\em Morita equivalent} if there is another local Kan simplicial manifold  $Z$ such that
both of the maps $ X\stackrel{\sim}{\leftarrow} Z
\stackrel{\sim}{\rightarrow} Y$ are hypercovers. In \cite[Section
  2]{z:tgpd-2}, we show that this definition does give an equivalence
relation. We call it {\em Morita equivalence} of local Kan simplicial manifolds.
\end{defi}

We also define {\em Morita
  equivalence} of Lie $n$-groupoids  
exactly in the same fashion  using hypercover of Lie $n$-groupoids.

Hypercover of Lie $n$-groupoids may also be  understood as a higher
analogue of pull-back of Lie groupoids. Let $X$ be a 2-groupoid and $Z_1\rra Z_0$ be two
manifolds with structure maps as in \eqref{eq:fd} up to
the level $n\leq 1$, and $f_n: Z_n \to X_n$ preserving the
structure maps $d^{n}_k$'s and $s^{n-1}_k$'s for $n\leq 1$. Then
$\hom(
\partial \Delta[n] , Z)$ still makes sense for $n\leq 1$. We further
suppose that
$f_0: Z_0 \twoheadrightarrow X_0$ (hence $Z_0\times Z_0
\times_{X_0\times X_0} X_1$ is a manifold) and $Z_1\thra Z_0\times
Z_0 \times_{X_0\times X_0} X_1$ are surjective submersions. That is to say  that the induced
map from $Z_k$ to the pull-back $\hom(
\partial \Delta[k] , Z)\times_{ \hom(\partial \Delta[k], X)}  X_k$
are surjective submersions for $k= 0, 1$. Then we form
\[Z_2= \hom(
\partial \Delta[2] , Z)\times_{ \hom(\partial \Delta[2], X)}  X_2,\]
which is a manifold (see \cite[Lemma 2.4]{z:tgpd-2}). 

Moreover there are $d^2_i: Z_2\to Z_1$ induced by the natural
projections $\hom(\partial \Delta[2] , Z)\to Z_1$; $s^1_i: Z_1 \to
Z_2$ by
\[ s^1_0(h)=(h,h,s^0_0(d^1_0(h)),s^1_0(f_1(h))), \quad
s^1_1(h)=(s^0_0(d^1_1(h)),h,h,s^1_1(f_1(h)));
\]
$m_i: \hom(\Lambda[3,i], Z) \to Z_2$ by for example
\[ m_0(( h_2, h_5, h_3, \bareta_1), (h_4, h_5, h_0, \bareta_2), (h_1, h_3, h_0, \bareta_3))=
( h_2, h_4, h_1, m_0(\bareta_1, \bareta_2, \bareta_3)), \] and
similarly for other $m$'s.
\[ \xymatrix{ &  & 0  & & \\
1 \ar[urr]^{h_0} & & & & 3 \ar[llll]^{h_4} \ar[llu]^{h_5} \ar[llld]^{h_2} \\
&  2 \ar[ruu]^{h_3} \ar[lu]^{h_1} & & & } \] Then $Z_2
\Rrightarrow Z_1 \rra Z_0$ is a Lie $2$-groupoid and we call it the
{\em pull-back 2-groupoid} by $f$. Moreover $f: Z\to X$ is an
equivalence with the natural projection $f_2: Z_2 \to X_2$.

\subsection{Lemmas}

What we wish to prove is:  If $X$ is already a Kan simplicial
manifold, then $X \stackrel{\sim}{\leftrightarrow}  Kan(X)$ are
Morita equivalent. It is very easy to prove for simplicial set. Since
the procedure of Kan replacement is basically to fill out horns, the
geometric realization of $Kan(X)$ and $X$ are homotopic to each
other. Since $X$ is Kan,  this is equivalent to \eqref{eqeq}.  However, the
missing tool of homotopy theory of  
simplicial manifolds (which do not form a model category, but building a certain machinery as a suitable replacement of model category
should be the eventual correct method to prove these lemmas.) 
prevents us to apply this proof directly. In fact, in the case of
simplicial sets,  one can easily obtain a morphism $\pi: Kan(X) \to X$
such that  the composition $X\to Kan(X) \to X$ is the identity. Then it is
straightforward to check that $Kan(X) \to X$ has the correct lifting
property. However, $\pi$ is not unique (basically it depends on the choice
of fillings in the Kan condition). Hence when generalized to a differential
category,   $\pi$ is in general not a continuous morphism. This forces us
to use another proof. Here we provide a proof
for Lie $2$-groupoids $W$. 

\begin{lemma}\label{lemma:kankan} If $W$ is a Lie $2$-groupoid, then
  $\tau_2(Kan(W))$ is a Lie 2-groupoid which is Morita equivalent to $W$. 
\end{lemma}
\begin{proof}

Usually, we do not have a direct map from  $\tau_2(Kan(W))$ to $W$
because  there is no (unique) multiplication map $W_1\times_{W_0}W_1 \to W_1$ (even when
there exists such a multiplication, we will encounter the issue of
surjective submersions). Hence we must construct a middle step. 

A more natural way  to describe this is to use the corresponding 
stacky groupoid  $\cG\rra W_0$,
where $\cG$ is
presented by the Lie groupoid $G_1 \rra G_0$, with $G_0= W_1$ and
$G_1 $ the set of bi-gons in $W_2$,  and the multiplication
$\cG\times_{W_0} \cG \to \cG$ is presented by bimodule
$E_m=W_2$. The bimodules of various compositions of multiplication from
various copies of $\cG$ to $\cG$ are presented by various fibre
product of $W_2$'s. For example, the bimodule $W_2 \times_{d_1, W_1, d_2} W_2$ with
the moment map $J_l$ to $W_1 \times_{W_0} W_1 \times_{W_0} W_1$ and $J_r$
to $W_1$, presents the multiplication
\[m\circ (m \times id):  (\cG \times_{W_0} \cG) \times_{W_0} \cG \to \cG. \]
To simplify the notation, we denote  a $k$-times fibre product as
$\square^{\times k}$ when it's clear from the context. We construct $ Z_0
= W_0 = W_0$ and, $Z_1$ is the disjoint union of these bibundles $W_2^{\times k}$ 
presenting different compositions of multiplication, 

\[
Z_1 = W_1 \sqcup W_2\sqcup (W_2
\sqcup W_2^{\times 2}  \sqcup  W_2^{\times 2}
 \sqcup  W_2^{\times 3} )
\sqcup \dots. \]
It is best to be understood as the following picture:\\
\psfrag{X1}{$Z_1$:} 
\centerline{\epsfig{file=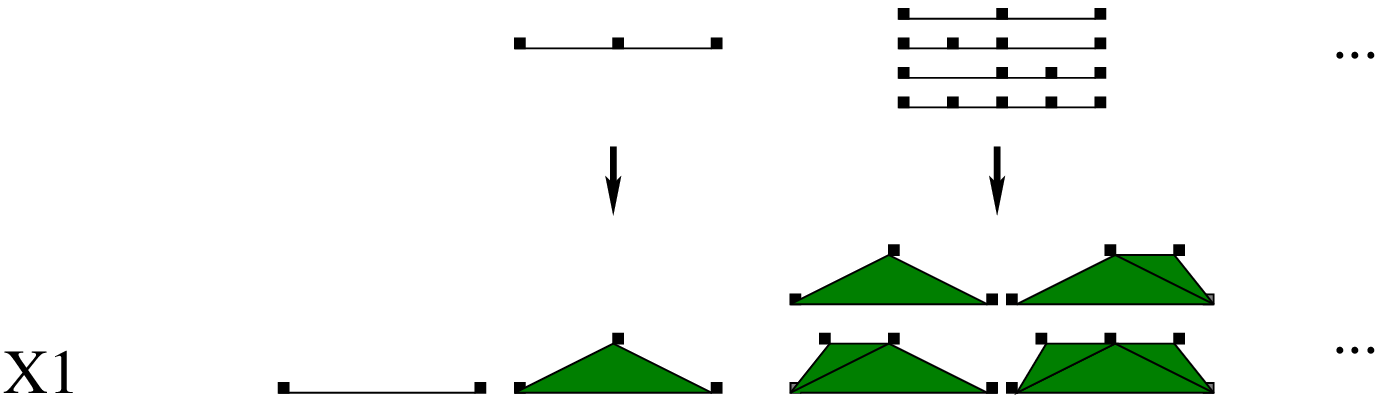, height=3cm, width=12cm}} 
 That is, we fill out horns in $Kan(W)_1$ by replacing
$W_1^{\times n}$ with $W_2^{\times{(n-1)}}$. The projections of $Z_1\to
Kan(W)_1$ and $Z_1 \to W_1$ are simply the disjoint union of the left and
right moment maps respectively. These projections are both surjective
submersions.

To show  that  $\tau_2(Kan(W))$ is Morita equivalent to $W$, we only
have to show that the pullback 2-groupoids on $Z$ are the same, that is 
\begin{equation}\label{eq:kanw-w}  Kan(W)_2 / \sim_2\times_{\hom(\pD[2],
 Kan( W))} \hom(\pD[2], Z)  \cong W_2 \times_{\hom(\pD[2],
  W)} \hom(\pD[2], Z). 
\end{equation}
If the map $p:M\to N$ is surjective and admit local section at any
point in $N$, then the pull-back groupoid $G_1 \times_M N \rra
G_0 \times_M N$ is free and proper if and only
the original groupoid $G_1 \rra G_0$ is so.  Since this is our case, the isomorphism \eqref{eq:kanw-w} automatically implies that
$Kan(W)/\sim_2$ is representable. By Prop. \ref{prop:tau2},
 $\tau_2(Kan(W))$ is a Lie 2-groupoid. 

We denote the two pullbacks by the map $Z_1\to W_1$ and
$Z_1 \to Kan(W)_1$ to $Z_1$ by $W|_Z$ and $Kan(W)|_Z$
respectively, and we construct morphisms
\[ \pi: (Kan(W)|_Z)_2 \to (W|_Z)_2, \quad \iota : (W|_Z)_2 \to (Kan(W)|_Z)_2 , \]
and prove $\pi \circ \iota = id$ and $\iota \circ \pi \sim id$ up to
something in $(Kan(W)|_Z)_3$. Then the above isomorphism follows
naturally. Notice that $Kan(W)$ is not a Lie 2-groupoid usually, but
pull-back described in Section \ref{sec:me} works also when $X$ is a
local Kan manifold. We form
$(X|_Z)_n = \hom ( sk_1(\Delta[n]) \to \Delta[n], Z\to X)$, where
$sk_1$ denotes of taking the 1-dimensional skeleton. By \cite[Lemma 2.4
]{z:tgpd-2}, $(X|_Z)_n$ are manifolds. Then it's easy to check that
$\tau_2(Kan(W)|_Z)=\tau_2(Kan(W))|_Z$. 

We first construct $\iota$.   Let $S$ be a simplicial
polygon with three marked points, namely a simplicial set constructed inductively
\[ \Delta[2]=S_0 \hookrightarrow S_1 \hookrightarrow S_2 \hookrightarrow
\dots \hookrightarrow S_i \dots,   \] 
by push-out $S_{i+1} = S_i \sqcup_{\Delta[1]} \Delta[2] $ and the three
marked points are the vertices of $S_0$. With these three marked
points, the $S_i$'s can be viewed as generalized triangles with their three sides
a concatenation  of line segments. In this sense, we also have the three facial maps $d^k$. 

We have a natural embedding
$W\hookrightarrow Kan(W)$, but this embedding does not give $W|_Z \to
Kan(W)|_Z$. In fact, take an element $ (w, \partial z) \in (W|_Z)_2 =
\sqcup_{i}\hom(S_i, W) $ for
a certain set of $S_i$'s, then $(w, \partial z) \notin (Kan(W)|_Z)_2$ since
$\partial w$ the boundary of $w$, is not $\partial z$ under the map
 $\partial_2 Z \to \partial_2 Kan(W)$. Here $\partial_k \square =
\hom(\pD[k], \square)$. To
construct $\iota$ we need to construct a morphism $\mu_i:\hom(S_i, W)
\to Kan(W)_2$ inductively, such that it commutes with the facial map
$d_k$ for $k=0,1,2$, 
\begin{equation}\label{eq:bd}\xymatrix{
\hom( S_i , W) \ar[r]^{d_k} \ar[d]^{\mu_i} & \hom(d^k S_i, W) \\
Kan(W)_2 \ar[ur]_{\partial}
} \end{equation} Then $\iota(w, \partial z):=(\mu_i(w), \partial z)$
where $w \in \hom(S_i, W)$.  

Step 1: We first prove the case $i=1$. We simplify the notation by $K:=Kan(W)$.
\begin{equation}\label{eq:fusion}\xymatrix{ \hom(S_1, W) = W_2
    \times_{d_k, W_1, d_1}  W_2 \ar[d]_{ \text{by}\; 
  W_1 \times_{W_0} W_1 \hookrightarrow K_2} \\
 W_2 \times_{d_k,
  W_1, d_1} W_2 \times_{W_1 \times_{W_0}  W_1} K_2
\ar[d]_{ \text{by}\;  W\hookrightarrow K }\\
\hom(\Lambda[3, k'], K) \ar[d]_{ Kan!(3, k')} \\
\hom(\Delta[3], K) \ar[d]_{d_{k'}} \\
K_2}\end{equation}
But this map does not commute with the facial map (see
\eqref{eq:bd}). For this purpose,  we only need to compose with the following one, 
\begin{equation}\label{eq:change}\xymatrix{K_2\to \big( K_2 \times_{W_1} (W_1\times_{W_0} W_1) \big) \times_{K_1
  \times_{W_0} W_1} W_1\times_{W_0} W_1  \ar[d]_{W_1 \times_{W_0} W_1
  \hookrightarrow K_2}\\
\hom(\Lambda[3, k''], K) \to K_2
}
\end{equation}
\psfrag{1}{{\tiny{$1$}}} \psfrag{0}{{\tiny{$0$}}} \psfrag{2}{{\tiny{$2$}}}
\psfrag{1'}{{\tiny{$1'$}}} \psfrag{S0}{{\tiny{$S_0$}}} \psfrag{K2}{{\tiny{$K_2$}}}
\psfrag{W11}{{\tiny{$W_1\times_{W_0}W_1$}}} \psfrag{good boundary}{good boundary}
\psfrag{eq1}{$\overset{\eqref{eq:fusion}}{\Rightarrow}$} \psfrag{eq3}{$\overset{\eqref{eq:change}}{\Rightarrow}$}
\centerline{\epsfig{file=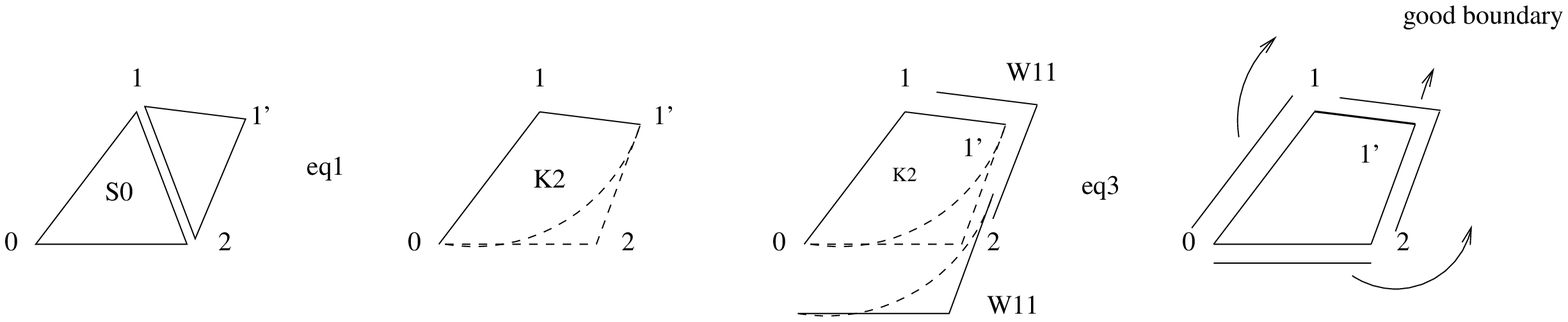, height=2.5cm} }

Step 2: Now suppose we have such a map $\hom(S_i ,W) \to K$, then
we can construct a map $\hom(S_{i+1}, W) \to K$ as below,
\[ \xymatrix{
\hom(S_{i+1} , W) = \hom(S_i, W) \times_{W_1} W_2 \ar[r]^{\partial}
\ar[d]_{\text{By Lemma \ref{lemma:wk}}  } &
\hom(\partial S_{i+1}, W) = \hom(\partial S_i, W)\times_{W_1}
\partial_2 W \ar[d] \\   
\hom(S_i, W) \times_{d_k, K_1, d_1} K_2 
\ar[r]^{\partial} \ar[d]_{\text{By} \; \hom(S_i, W) \to K_2 } & \hom(\partial S_i, W)
\times_{K_1} \partial_2 K \ar[d]\\
K_2\times_{d_k, K_1, d_1} K_2 =\hom(S_1, K) 
\ar[d]_{\text{Similarly as Step 1,  replace $W$ by $K$}}
\ar[r]^{\partial} & \hom(\partial S_1, K)
\\K_2 \ar[ur]_{\partial} 
}
\]
\psfrag{0}{$0$}\psfrag{1}{$1$}\psfrag{2}{$2$} \psfrag{eq1}{$\overset{\text{Lemma
  \ref{lemma:wk}}}{\Rightarrow}$} \psfrag{eq2}{$\overset{\text{Viewed as}}{\Rightarrow}$}
\psfrag{eq3}{$\overset{\text{As Step 1}}{\Rightarrow}$}
\centerline{\epsfig{file=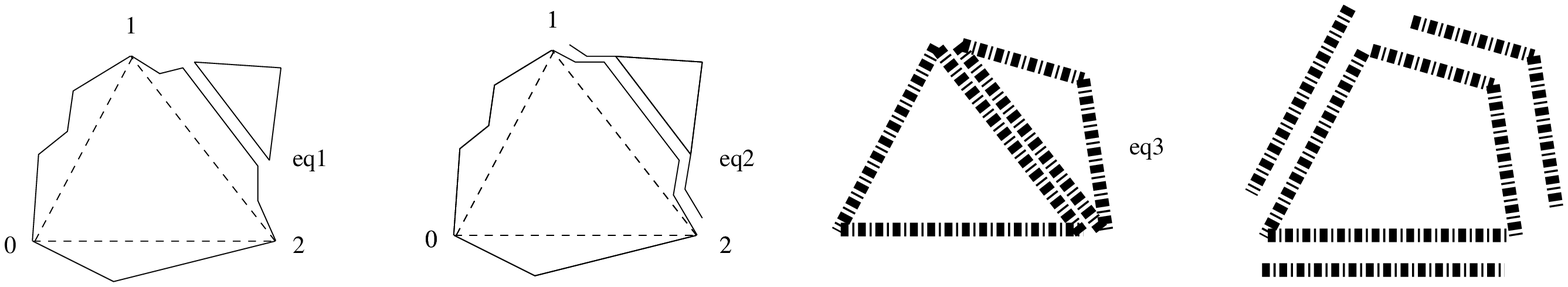, height=2.5cm}}

\begin{lemma}\label{lemma:wk} There is a natural morphism $K_1 \times_{W_1}
  W_2 \to K_2$.  
\end{lemma}
We leave this to the readers as an exercise involving the Kan condition of $K$. 

To construct $\pi$, we first construct a local morphism $f: Kan(W)
\to W  $
 inductively.  
The first step is to construct $f^1_n$ by  
\[ 
\begin{split}
W^{1}_n 
& = \hom\left( \D[n], W \coprod_{\coprod \L[k,j] \times
  \hom(\L[k,j], W)} \D[k] \times \hom(\L[k,j], W) \right) \\
& \to W_n \coprod_{\dots} \hom(\D[n], \D[k]) \times \hom(\D[k], W) \\
& \to W_n .
\end{split}
\] 
In the second last step we use the strict Kan condition
$\hom(\L[k,j], W) \cong  W_k$ when $k\ge 2$ and we choose a
local section $\hom(\L[2,1], W) \to  W_2$ when $k=2$. The last step
follows from the composition $\hom(\D[n], \D[k]) \times \hom(\D[k], W)
\to W_n$ and thus both spaces in the  coproduct have a natural map to $W_n$.

Suppose that $f^\beta : W^\beta 
\to W$ is constructed. Then $f^{\beta+1}$ is the composition of the
following natural morphisms
\[
\begin{split}
W^{\beta+1}_n &=\hom \left( \D[n], W^\beta \coprod_{\coprod \L[k,j] \times
  \hom(\L[k,j], W^\beta)} \D[k] \times \hom(\L[k,j], W^\beta) \right)
\\
& \to \hom\left( \D[n], W \coprod_{\coprod \L[k,j] \times
  \hom(\L[k,j], W)} \D[k] \times \hom(\L[k,j], W) \right) \\
& = W^1_n \xrightarrow{f^1_n} W_n .
\end{split}
\]
Then $f$ is the colimit of $f^\beta$. 

More geometrically, if we view  an element in
$Kan(W)_2$ as a set of  small
triangles of $W_2$ touching together,  $f_2$ is basically
to compose these small triangles into a big one in $W_2$ with a choice
of  filling for each $W_1\times_{W_0} W_1$, which is given by $f^1$.  \\
\psfrag{in Y}{$f(x)$} \psfrag{trivial filling}{trivial
  filling} \psfrag{element in Z}{element in $Z_1$}
\psfrag{yz}{$f(x) \circ z$}
\centerline{\epsfig{file=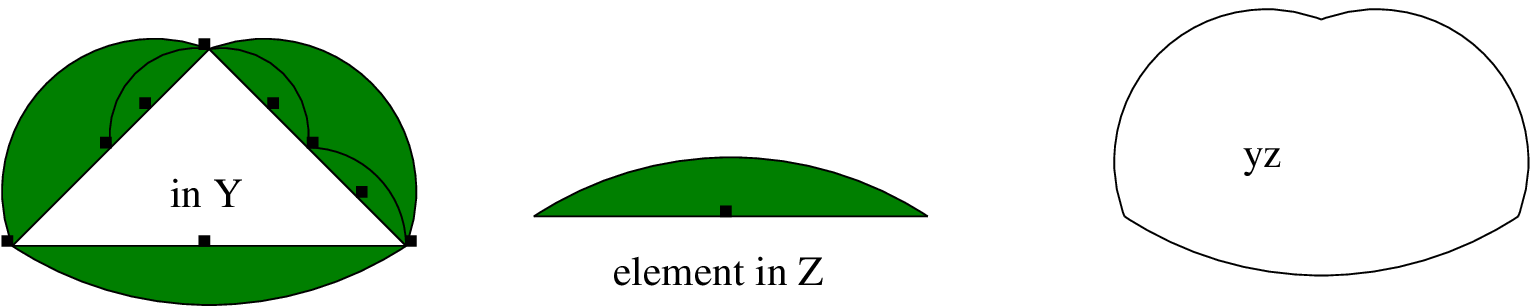, height=3cm, width=12cm} }
Now when we make a choice of fillings for a $W_1\times_{W_0}W_1$ on
the boundary, instead of choosing some filling given by $f^1$, we
choose the element $\partial z \in \hom(\pD[2], Z) $, then this
element in $W_2$ is denoted by $f_2(x) \circ \partial z$.  

Thus $f_2$ induces a map
\begin{equation}\label{eq:varphi} Kan(W)_2  \times_{\hom(\pD[2],
 Kan( W))} \hom(\pD[2], Z)  \xrightarrow{\pi} W_2 \times_{\hom(\pD[2],
  W)} \hom(\pD[2], Z). 
\end{equation}
as $(x, \partial z ) \mapsto (f_2(x) \circ \partial z, \partial z)$.
In
Lemma \ref{lemma:tree},
we give a combinatorial proof that this map does not depend on the
choice of fillings. Hence we obtain a well-defined global map
$\pi$.

Then it is not hard to see that $\pi \circ \iota = id$ since $\pi$ is
exactly the opposite procedure of $\iota$. 

The procedure to form $\pi$ and $\iota$ is basically to use $Kan!(3,
j)$ to compose (for example \eqref{eq:fusion}), hence $\iota \circ \pi$ and $ id$ differ by something
in $Kan(W)_3$.  
\end{proof}

\begin{lemma}\label{lemma:tree}
The map $\pi$ does not depend on the choice of sections in  the construction of $f$. 
\end{lemma}
\begin{proof}
We denote an element in $\eta \in Kan(W)_2 $ by a bicolored tree\\
\centerline{\epsfig{file=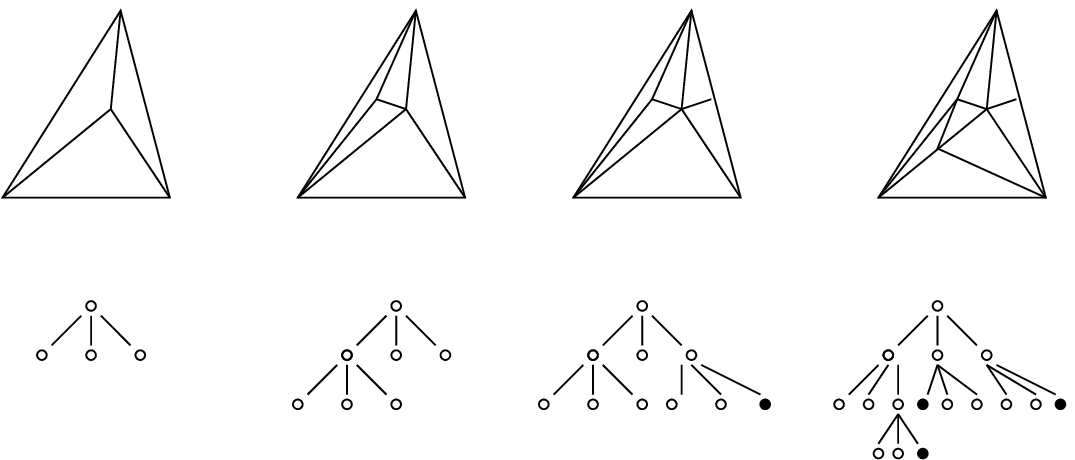, height=3.8cm, width=12cm} }
A point is black if it represents a weird triangle, i.e. a triangle
comes from the first copy of $X_1 \times_{X_0} X_1$ which serves as an
artificial filling; otherwise, it is
white. For the other two copies of $X_1\times_{X_0}X_1$, they are
degenerate ones and can not glue directly with a normal triangle in
$W_2$. Since degenerate elements play the role of identities in
composition $\pi$, we here ignore them.  We prove the result by induction on the number of generations
and the number of points in the youngest generation. It is obvious
for the initial case. 

Now take three siblings points in the youngest
generation, if all of them are white, then we use $Kan(3, j)$
without a choice and we end up with an element $\eta' \in Kan(W)_2$
which has a fewer number of generations
or a fewer number of points in the youngest generation. Done! 

If one of the three siblings is black, then there is precisely one
black one in these three siblings, which we denote by $x$. Since $\hom(\pD[2], Z)$ will give the fillings for
the weird triangles on the border of $\eta$, to show the independence,
we only have to deal with the inner triangles. Then some
ancestor of $x$ must have a black descendant $y$, 
because a weird triangle must lie
on the side of another triangle (which is the parent of $y$). 

The simplest situation in this case 
is when the other black descendant is a (true) cousin (namely their
direct ancestors are siblings) \\
\begin{equation}\label{pic:simple}
\psfrag{x}{$x$}\psfrag{y}{$y$}
\centerline{\epsfig{file=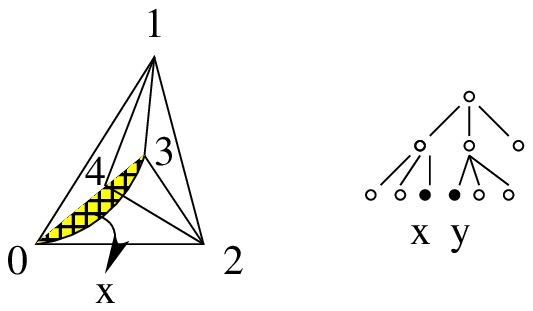, height=3cm, width=10cm} }
\end{equation}
We choose a triangle $\eta_{034}$, and we are given triangle $\eta_{014}
$, $\eta_{134}
$, $\eta_{123}
$, $\eta_{234}
$,  $\eta_{024}
$. The procedure is to compose $\eta_{034}$,  $\eta_{014}
$, $\eta_{134}
$ first to obtain $\eta_{013}$ by $Kan(3, 3)$; then to compose $\eta_{034}$,  $\eta_{024}
$, $\eta_{234}
$ secondly to obtain $\eta_{023}$ by $Kan(3,3)$; finally to compose $\eta_{013}$,
$\eta_{023}$, and $\eta_{123} $ to obtain $\eta_{012}$ by $Kan(3,
3)$. These can be viewed as
multiplications for 2-groupoid \cite[Section 2.3]{z:tgpd-2}. By
associativity of such multiplications (or equivalent $Kan(3, j)!$ and
$Kan(4, j)!$), we can obtain the same
$\eta_{012}$ by another order of composition, namely we use $Kan(3,2)$
first to obtain $\eta_{124}$ then $Kan(3, 3)$.  Since the second way to
compose does not depend on the choice of $\eta_{034}$, our final
result $\eta_{012}$ does not depend on the choice either.  Hence by
the induction hypothesis, we will choice-independently end up with an
element $f_2(x)\circ \partial z$.  

We might meet more complicated situations, namely the other black
descendant $y$ is a more remote cousin,  but we can reduce them to
the simple situation above: \\
\psfrag{x}{\tiny{$x, \xi$}}\psfrag{y}{\tiny{$y, \xi$}}\psfrag{xx}{$x', \xi'$}\psfrag{yy}{$y', \xi''$}\psfrag{z}{$ \zeta$}
\centerline{\epsfig{file=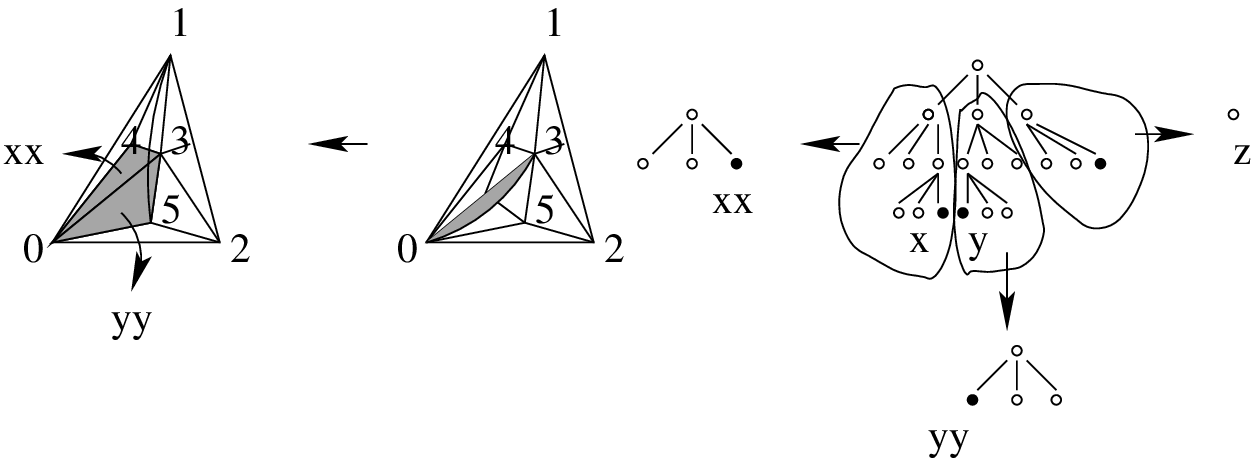, height=4.5cm, width=12cm} }
We choose a triangle $\xi\in W_2$ corresponding to $x$ in the
youngest generation. Then we use $Kan(3, j)$ to compose to obtain
$\xi'$ and  $\xi''$ corresponding to $x'$ and $y'$ respectively. As
shown in the picture, the true cousins $x', y'$ are ancestors of $x$
and $y$ respectively. We
also do the same for the other branch and obtain $\zeta$. 
During this procedure, we might
have to make other choices of fillings for other black points. But it
does not matter,  since our point is to show the independence on the choice $\xi$.  Then we are
again in a similar situation as of \eqref{pic:simple}.   We choose
a filling $\eta_{234}$. By $Kan(4, j)!$ the final result does not
depend on the order of composition.  Then by $Kan(3, 2)$ 
we first obtain $\eta_{124}$; by the induction hypothesis, we obtain $\eta_{024}$ independent of the
choice of $\xi$ because $\eta_{024}$ contains fewer descendants;
finally, by $Kan(3, 3)$ we obtain $\eta_{012}$, which is independent
of the choice of filling $\xi$. 

\end{proof}

\begin{lemma}
If both $X$ and $Y$ are Lie $n$-groupoids, then a hypercover of 
local Kan simplicial manifolds $X \lhc Y$ is automatically a hypercover of
Lie $n$-groupoids.
\end{lemma}
\begin{proof}
Since $X$ is a Lie $n$-groupoid, we have the composed morphisms $\hom(\L[n+1,
j], X) \cong X_{n+1} \xrightarrow{d_j} X_n \to Y_n$ and $\hom(\L[n+1,
j], X) \cong X_{n+1} \xrightarrow{d_j} X_{n} \xrightarrow{\partial}
\hom(\pD[n], X)$. This gives us a map $\hom(\L[n+1, j], X) \xrightarrow{p}
\hom(\pD[n]\to \D[n], X \to Y)$. With this map, we rewrite
\[\hom(\pD[n+1]\to \D[n+1], X \to Y) \cong \hom(\L[n+1, j], X)
\times_{\hom(\pD[n]\to \D[n], X \to Y)} X_n. \] 
Since $X\lhc Y$ as local Kan simplicial manifolds, the following map 
\[ X_{n+1}\cong   \hom(\L[n+1, j], X) \to \hom(\L[n+1, j], X)
\times_{\hom(\pD[n]\to \D[n], X \to Y)} X_n \]
is a surjective submersion. This implies that $X_n \to \hom(\pD[n]\to \D[n], X \to Y)$ is 
injective. However, $X_n \to \hom(\pD[n]\to \D[n], X \to Y) $ is a surjective
submersion by the condition of hypercovers. Hence $X_n \cong
\hom(\pD[n]\to \D[n], X \to Y)  $, which shows $X\lhc Y$ as Lie $n$-groupoids. 
\end{proof}

This implies
\begin{cor}
Two Lie $n$-groupoids $X\me Y$ are Morita equivalent as local  Kan simplicial
manifolds if and only if they are Morita equivalent as Lie $n$-groupoids.
\end{cor}

\begin{lemma}\label{lemma:trun2} If $\phi: X \lhc Y $ is a hypercover of
local Kan  simplicial manifolds, and if $Kan(X)_2/\sim_2$ is
representable,  then both $\tau_2(Kan(X))$ and $\tau_2(Kan(Y))$ are
Lie 2-groupoids and the induced map $\tau_2(Kan(X))
\to \tau_2(Kan(Y))$ is  a hypercover of Lie 2-groupoids. 
\end{lemma}
\begin{proof}
We first show that if $\phi: K \to K'$ is a hypercover of Kan simplicial
sets (i.e. \eqref{eq:pb} is surjective instead of a surjective submersion), then the natural map
\begin{equation}\label{eq:kanx-kany}
\tau_n(K)_n \xrightarrow{f} \hom(\pD[n]\to \D[n], \tau_n(K) \to \tau_n(K')),  
\end{equation}
is an isomorphism. 
Notice that the right hand side is simply $\hom(\pD[n],
K)\times_{\hom(\pD[n], K')} \hom(\D[n], \tau_n(K'))  $. Thus we have a commutative diagram
\[
\xymatrix@C=3.5cm{K_n \ar@{->>}[r] \ar@{->>}[d]_g &\hom(\pD[n]\to \D[n], K \to K')\ar@{->>}[d] \\
\tau_n(K)_n \ar[r]^f & \hom(\pD[n]\to \D[n], \tau_n(K) \to \tau_n(K')),}
\]
where $\thra$ denotes surjective maps. Then $f$ must be  surjective because $f\circ g$ 
being  surjective  implies $f$ being  surjective.

Take $(\delta x_n, [y_n]) \in \hom(\pD[n],
K)\times_{\hom(\pD[n], K')} \hom(\D[n], \tau_n(K')) $. If both
$[x_n], [x'_n]\in \tau_n(K)_n$ map to $(\delta x_n, [y_n])$, that is
the boundary $\partial
x_n = \partial x'_n = \delta
x_n \in \hom(\pD[n], K)$ and $\phi_n(x'_n)=y'_n \sim y_n \sim y''_n
=\phi_n(x_n)$, then 
$y''_n$ and $y'_n$ differ by a certain element $y_{n+1}\in K'_{n+1}$. Since 
\[ K_{n+1}  \thra \hom(\pD[n+1]\to \D[n+1], K \to K'),   \] 
is  surjective, there exists $x_{n+1}$ such that
$\phi_{n+1}(x_{n+1})=y_{n+1}$ and  $x_n,
x'_n$ differ by $x_{n+1}$. This proves that $[x_n]=[x'_n]\in
\tau_n(K)_n$. Hence $f$ is also injective.

If the map $p:M\to N$ is surjective and admits local section at any
point in $N$, then the pull-back groupoid $G_1 \times_M N \rra
G_0 \times_M N$ is free and proper if and only
the original groupoid $G_1 \rra G_0$ is so.  Since this is our case,
the isomorphism \eqref{eq:kanx-kany}, when applied to $n=2$ and
$K=Kan(X)$ $K'=Kan(Y)$,  implies that
$Kan(Y)/\sim_2$ is representable. Hence $\tau_2(Kan(X))$ and $\tau_2(Kan(Y))$
are Lie 2-groupoids by Prop. \ref{prop:tau2}.

Now we only need to verify that the morphism
\[
\tau_2(Kan(X))_m  \to \hom(\pD[m]\to \D[m], \tau_2(Kan(X))\to \tau_2( Kan(Y)))
\]
is a surjective submersion for $m=0, 1$. For $m=0$ it is implied by
$X_0 \to Y_0$ being a surjective submersion. For $m=1$, by induction, we need to show that the natural map
\begin{equation}\label{eq:aim} X^{\beta+1}_{1} \to \hom(\pD[1]\to \D[1], X^{\beta+1} \to
Y^{\beta+1} ), \end{equation}
is a surjective submersion supposing the same is true for $\beta$. We have
\[  
\begin{split}
X^{\beta+1}_{1}&=X^\beta_1  \sqcup \hom(\L[2,1],X^\beta), \quad Y^{\beta+1}_{1}=Y^\beta_1  \sqcup\hom(\L[2,1],
Y^\beta)
.
\end{split}
 \]

The right hand side of \eqref{eq:aim} decomposes
into two terms $I, II$ according to the decomposition of
$Y^{\beta+1}_1$,   
\[
\begin{split}
I&=\hom(\pD[1], X^\beta) \times_{\hom(\pD[1], Y^\beta)} \hom(\D[1],
Y^\beta) \\
II&=\hom(\pD[1], X^\beta) \times_{\hom(\pD[1], Y^\beta)} \hom(\L[2,
  1], Y^\beta) 
\end{split}
\] 
By the induction hypothesis, $X^\beta_1 \to I
$ is a surjective submersion. Further by \cite[Lemma 2.5]{z:tgpd-2} (take
$S=T=\L[2, 1]$, and $T'=\partial \D[1]$),  
\[ \hom(\L[2, 1], X^\beta) \to II\]
is a surjective submersion.
Thus \eqref{eq:aim} is a surjective submersion.
\end{proof}

\bibliographystyle{habbrv}
\bibliography{../../../bib/bibz.bib}

\def\cprime{$'$} \def\cprime{$'$} \def\cprime{$'$} \def\cprime{$'$}
\begin{thebibliography}{10}

\bibitem{SGA4}
{\em Th\'eorie des topos et cohomologie \'etale des sch\'emas. {T}ome 1:
  {T}h\'eorie des topos}.
\newblock Lecture Notes in Mathematics, Vol. 269. Springer-Verlag, Berlin,
  1972.
\newblock S{\'e}minaire de G{\'e}om{\'e}trie Alg{\'e}brique du Bois-Marie
  1963--1964 (SGA 4), Dirig{\'e} par M. Artin, A. Grothendieck, et J. L.
  Verdier. Avec la collaboration de N. Bourbaki, P. Deligne et B. Saint-Donat.

\bibitem{cafe}
A.~S. Cattaneo and G.~Felder.
\newblock Poisson sigma models and symplectic groupoids.
\newblock In {\em Quantization of singular symplectic quotients}, volume 198 of
  {\em Progr. Math.}, pages 61--93. Birkh\"auser, Basel, 2001.

\bibitem{cf}
M.~Crainic and R.~L. Fernandes.
\newblock Integrability of {L}ie brackets.
\newblock {\em Ann. of Math. (2)}, 157(2):575--620, 2003.

\bibitem{friedlander}
E.~M. Friedlander.
\newblock {\em \'{E}tale homotopy of simplicial schemes}, volume 104 of {\em
  Annals of Mathematics Studies}.
\newblock Princeton University Press, Princeton, N.J., 1982.

\bibitem{getzler}
E.~Getzler.
\newblock {Lie theory for nilpotent L-infinity algebras},
  arxiv:math.AT/0404003.

\bibitem{henriques}
A.~Henriques.
\newblock Integrating {$L\sb \infty$}-algebras.
\newblock {\em Compos. Math.}, 144(4):1017--1045, 2008.

\bibitem{lurie}
J.~Lurie.
\newblock {Higher Topos Theory}, arXiv:math/0608040v4 [math.CT].

\bibitem{may}
J.~P. May.
\newblock {\em Simplicial objects in algebraic topology}.
\newblock Chicago Lectures in Mathematics. University of Chicago Press,
  Chicago, IL, 1992.
\newblock Reprint of the 1967 original.

\bibitem{hs}
J.~Mr{\v{c}}un.
\newblock {\em Stablility and invariants of {H}ilsum-{S}kandalis maps}.
\newblock Dissertation, Utrecht University, Utrecht, 1996.

\bibitem{quillen:ha}
D.~G. Quillen.
\newblock {\em Homotopical algebra}.
\newblock Lecture Notes in Mathematics, No. 43. Springer-Verlag, Berlin, 1967.

\bibitem{tz}
H.-H. Tseng and C.~Zhu.
\newblock Integrating {L}ie algebroids via stacks.
\newblock {\em Compos. Math.}, 142(1):251--270, 2006.

\bibitem{van-est:local}
W.~T. van Est and M.~A.~M. van~der Lee.
\newblock Enlargeability of local groups according to {M}al\cprime cev and
  {C}artan-{S}mith.
\newblock In {\em Action hamiltoniennes de groupes. {T}roisi\`eme th\'eor\`eme
  de {L}ie ({L}yon, 1986)}, volume~27 of {\em Travaux en Cours}, pages 97--127.
  Hermann, Paris, 1988.

\bibitem{severa:diff}
P.~\v{S}evera.
\newblock {L infinity algebras as 1-jets of simplicial manifolds (and a bit
  beyond) }, arXiv:math/0612349 [math.DG].

\bibitem{z:tgpd-2}
C.~Zhu.
\newblock { n-groupoids and stacky groupoids},
  arxiv:math.DG/0801.2057[math.DG].

\bibitem{z:lie2}
C.~Zhu.
\newblock {Lie II theorem for Lie algebroids via stacky Lie groupoids},
  arxiv:math/0701024 [math.DG].

\bibitem{z:tgpd}
C.~Zhu.
\newblock {Lie n-groupoids and stacky Lie groupoids}, arxiv:math.DG/0609420.

\end{thebibliography}

\end{document}